\documentclass[12pt,reqno]{amsart}
\usepackage{enumerate, latexsym, amsmath, amsfonts, amssymb, amsthm, graphicx, color}
\usepackage{a4}
\usepackage{amsmath}
\usepackage{amsfonts}
\usepackage{amsthm}
\usepackage{amssymb}
\usepackage{graphicx}
\def\pmod #1{\ ({\rm{mod}}\ #1)}
\def\Z{\Bbb Z}
\def\N{\Bbb N}

\def\l{\left}
\def\r{\right}
\def\bg{\bigg}
\def\({\bg(}
\def\){\bg)}
\def\t{\text}
\def\f{\frac}

\def\mo{{\rm{mod}\ }}

\def\ls{\leqslant}
\def\gs{\geqslant}

\def\sm{\setminus}
\def\bi{\binom}

\def\eq{\equiv}

\def\Proof{\noindent{\it Proof}}

\makeatletter
\@namedef{subjclassname@2020}{%
  \textup{2020} Mathematics Subject Classification}

 \makeatletter

\newtheorem{theorem}{Theorem}[section]

\newtheorem{lemma}[theorem]{Lemma}

\theoremstyle{definition}

\newtheorem{remark}{Remark}[section]

\def \tmz {\begin{itemize}}
\def \etmz {\end{itemize}}

\numberwithin{equation}{section}
\begin{document}
\hbox{Polished version for submission}
\medskip

\title[Universal sums  via products of Ramanujan's theta functions]
{Universal sums via products of \\ Ramanujan's theta functions}

\author[N. A. S. Bulkhali and  Z.-W. Sun]{Nasser Abdo Saeed Bulkhali and Zhi-Wei Sun}
	
\address{(Nasser Abdo Saeed Bulkhali) Department of Mathematics, Faculty of Education-Zabid, Hodeidah University, Hodeidah, Yemen}
\email{nassbull@hotmail.com}

\address {(Zhi-Wei Sun) School of Mathematics, Nanjing University,
	Nanjing 210093, People's Republic of China}
\email{zwsun@nju.edu.cn}
\date{}

\keywords{Ramanujan's theta function,  sums of squares, triangular numbers, ternary quadratic forms.
\newline \indent The initial version of this paper was posted to arXiv as a preprint with the ID
{\tt arXiv:2410.14605}. The second author is supported by the Natural Science Foundation of China (grant 12371004).}

\begin{abstract}
An integer-valued polynomial $P(x,y,z)$ is said to be universal (over $\Z$) if each nonnegative integer can be written as $P(x,y,z)$ with $x,y,z\in\Z$. In this paper, we mainly introduce a new technique to determine the universality of some sums in the form $x(a_1x+a_2)/2+y(b_1y+b_2)/2+z(c_1z+c_2)/2$
(with $a_1-a_2,b_1-b_2,c_1-c_2$ all even)
conjectured by Sun,  using various identities of Ramanujan's theta functions.
For example, we prove that $x(3x+1)+y(3y+2)+2z(3z+2)$
and $x(4x+r)+y(3y+1)/2+z(7z+1)/2\ (r=1,3)$ are universal.
\end{abstract}

\maketitle

\section{Introduction}\label{Introduction}

Let $\N=\{0,1,2,\ldots\}$ and $\Z^+=\{1,2,3,\ldots\}$.
For an integer-valued polynomial $P(x_1,\ldots,x_n)$, we call it {\it universal} (over $\Z$) if
$$\{P(x_1,\ldots,x_n):\ x_1,\ldots,x_n\in\Z\}=\N,$$
and call $P(x_1,\ldots,x_n)$ {\it universal over $\N$} if
$$\{P(x_1,\ldots,x_n):\ x_1,\ldots,x_n\in\N\}=\N.$$
When
$$\N\sm\{P(x_1,\ldots,x_n):\ x_1,\ldots,x_n\in\Z\}$$
is finite, we say that $P(x_1,\ldots,x_n)$ is {\it almost universal} (over $\Z$).

For each integer $m\gs3$, those numbers
$$p_m(n)=(m-2)\bi n2+n=\f{(m-2)n^2-(m-4)n}2\ (n=0,1,2,\ldots)$$
are called {\it $m$-gonal numbers}, and those numbers $p_m(x)$ with $x\in\Z$
are usually called {\it generalized $m$-gonal numbers}. Note that
\begin{equation}\label{p-}p_m(-x)=\f{x((m-2)x+m-4)}2.
\end{equation}
In particular,
\begin{gather*} p_3(x)=\f{x(x+1)}2,\ p_3(-x)=\f{x(x-1)}2,\ p_4(x)= p_4(-x)=x^2,\ \\p_5(x)=\f{x(3x-1)}2,
\ p_5(-x)=\f{x(3x+1)}2,
\\p_6(x)=x(2x-1)=p_3(2x-1),\ p_6(-x)=x(2x+1)=p_3(2x),
\\ p_7(x)=\f{x(5x-3)}2,\ p_7(-x)=\f{x(5x+3)}2,\ p_8(x)=x(3x-2),\ p_8(-x)=x(3x+2).
\end{gather*}
Obviously,
\begin{equation}\label{p6}\{p_6(x):\ x\in\Z\}=\{x(2x+1):\ x\in\Z\}=\{p_3(x):\ x\in\Z\}=\{p_3(n):\ n\in\N\}.
\end{equation}
Fermat claimed that $p_m(x_1)+\cdots+p_m(x_m)$ is universal over $\N$. This was confirmed by
Lagrange, Gauss and Cauchy for the cases $m=4$, $m=3$, and $m\gs5$, respectively.
In 1862 Liouville (cf. \cite[p. 23]{Dickson}) determined all universal sums $ap_3(x)+bp_3(y)+cp_3(z)$ with $a,b,c\in\Z^+$,
namely, those triples $(a,b,c)$ with $a\ls b\ls c$ such that $ap_3(x)+bp_3(y)+cp_3(z)$ is universal are
$$(1,1,1),\ (1,1,2),\ (1,1,4),\ (1,1,5),\ (1,2,2),\ (1,2,3),\ (1,2,4).$$
 In $2007$, Z.-W. Sun \cite{Sun} initiated the study of universal sums of the types
 $ax^2+by^2+cp_3(z)$ and $ax^2+bp_3(y)+cp_3(z)$ with $a,b,c\in\Z^+$, and the investigation was finished
 via the three papers \cite{Sun,GPS,OS}.

 In 2010 B. Kane and Sun \cite{KaneSun} used modular forms and the theory of ternary quadratic forms to characterize those tuples $(a,b,c)\in(\Z^+)^3$
 such that $ax^2+by^2+cp_3(z)$ is almost universal; they also investigated those tuples $(a,b,c)\in(\Z^+)^3$  such that $ax^2+bp_3(y)+cp_3(z)$ is almost universal, as well as those tuples $(a,b,c)\in(\Z^+)^3$  such that $ap_3(x)+bp_3(y)+cp_3(z)$ is almost universal.
 Via solving a conjecture of Kane and Sun~\cite{KaneSun}, W. K. Chan and B.-K. Oh \cite{ChanOh} gave a complete characterization of those triples $(a,b,c)\in(\Z^+)^3$ with $ax^2+bp_3(y)+cp_3(z)$
 (or $ap_3(x)+bp_3(y)+cp_3(z)$) almost universal.

 In 2015 Sun \cite{Sun2015} investigated universal sums of the type $ap_i(x)+bp_j(y)+cp_k(z)$ with $a,b,c\in\Z^+$ and $i,j,k\in\{3,4,5,\ldots\}$. In particular,  for $k\in\{4,5,7,8,9,\ldots\}$ and $a,b,c\in\Z^+$ with $a\ls b\ls c$, Sun \cite{Sun2015} showed that if $ap_k(x)+bp_k(y)+cp_k(z)$ is universal then $k=5$, $a=1$  and $(b,c)$ is among the following $20$ ordered pairs:
 $$(1,c)\, (c\in\{1,\ldots,10\}\sm\{7\}),\ (2,c)\, (c\in\{2,3,4,6,8\}),
(3,c)\, (c=3,4,6,7,8,9).$$
For these ordered pairs $(b,c)$, it is now known (cf. \cite{Guy,Sun2015,O}) that $p_5(x)+bp_5(y)+cp_5(z)$ is indeed universal.

In $2017$,  Sun \cite{Sun1} proved that, for  integers $0 \ls b \ls c \ls d \ls a$ with $a >2$, the sum $x(ax +b) +y(ay+c) +z(az+d)$ is universal if and only if the quadruple $(a, b, c, d)$ is among
$$(3, 0, 1, 2),~ (3, 1, 1, 2),~ (3, 1, 2, 2),~ (3, 1, 2, 3),~ (4, 1, 2, 3).$$
Sun \cite{Sun1} also proved that $x(ax+1)+y(by +1)+z(cz +1)$ is universal if $(a, b, c)$ is among the following seven triples:
$$(1, 2, 3),~ (1, 2, 4),~ (1, 2, 5),~ (2, 2, 4),~ (2, 2, 5),~ (2, 3, 3),~ (2, 3, 4).$$
J. Ju and Oh \cite{JuOh} showed that $x(ax + 1) + y(by + 1) + z(cz + 1)$ is also
universal for $(a, b, c) = (2, 2, 6),~ (2, 3, 5),~ (2, 3, 7)$ as conjectured
by Sun \cite{Sun1}.

For convenience, for any $a_1,b_1,c_1\in\Z^+$ and $a_2,b_2,c_2\in\Z$ with $a_1-a_2,b_1-b_2,c_1-c_2$
 all even, whenever
 \begin{equation}\label{xyz}\f{x(a_1x+a_2)}2+\f{y(b_1y+b_2)}2+\f{z(c_1z+c_2)}2
\end{equation}
is universal (resp., almost universal), we call the tuple $(a_1,a_2,b_1,b_2,c_1,c_2)$  universal
(resp., almost universal).
For any $a\in\Z^+$, by \eqref{p6} we have
$$\{ap_3(x):\ x\in\Z\}=\l\{\f{x(ax+a)}2:\ x\in\Z\r\}=\left\{\f{x(4ax+2a)}2:\ x\in\Z\right\}.$$
Sun \cite{Sun2,Sun3} investigated those {\it standard} universal tuples $(a_1,a_2,b_1,b_2,c_1,c_2)\in\N^6$
with $a_1-a_2, b_1-b_2,c_1-c_2\in2\Z^+$, $a_1\gs b_1\gs c_1\gs2$,
and $a_2\gs b_2$ if $a_1=b_1$, and $b_2\gs c_2$ if $b_1=c_1$; he
showed that such tuples $(a_1,a_2,b_1,b_2,c_1,c_2)$ are among the $12082$ tuples listed in \cite[Appendix]{Sun2}.
 Those tuples in the list marked with $*$ have been proved to be universal,
 and Sun \cite{Sun2,Sun3} conjectured that all the remaining tuples in the list are indeed universal.
 In this direction, H.-L. Wu and Sun \cite{WuSun} made certain progress.
 As Sun \cite{Sun1} proved that $x(2x+1)+y(3y+1)+z(4z+1)$ is universal, the tuple $(8,2,6,2,4,2)$
 not marked with $*$ in the list is indeed universal. Note that all previous works in this direction
 depend heavily on the theory of ternary quadratic forms.

 In Section 2, we will confirm that some open tuples listed in \cite[Appendix]{Sun2} are indeed universal
 by using the theory of ternary quadratic forms.

Ramanujan's general theta function is defined by
\begin{equation}\label{2.1}f(a,b):= \sum_{n=-\infty}^{\infty}{{a}^{n(n+1)/2}}{{b}^{n(n-1)/2}}\ \ \ \t{with}\ |ab|<1.\end{equation}
In particular, Ramanujan \cite[Entry 22]{Adiga3} introduced $\varphi(q)$ and $\psi(q)$  for $|q|<1$:
\begin{align}\label{2.7}\varphi(q)&:= f(q,q) =  \sum_{n=-\infty}^{\infty} q^{n^2},\\
\label{2.8}\psi(q)&:= f(q,q^3) =  \sum_{n=-\infty}^{\infty} {q^{n(2n+1)}}=  \sum_{n=0}^{\infty} {q^{n(n+1)/2}}.
\end{align}
For various formulas of Ramanujan's theta functions, one may consult the excellent introductory book
of B.C. Berndt \cite[Chapter 1]{Berndt}.

In this paper we find that many universal sums can be studied via identities of Ramujan's theta functions.
We establish the following novel theorem connecting universal sums with Ramujan's theta functions.

\begin{theorem}\label{Th1.1} Let $i,j,s,t,u,v\in\Z$ with $i+j,s+t,u+v$ positive.
Suppose that
\begin{equation}\label{=}f(q^i,q^j)f(q^s,q^t)f(q^u,q^v)=\sum_{r=1}^k m_rq^{r-1}f(q^{ka_{1r}},q^{ka_{2r}})
f(q^{kb_{1r}},q^{kb_{2r}})f(q^{kc_{1r}},q^{kc_{2r}}),
\end{equation}
where $a_{1r},a_{2r},b_{1r},b_{2r},c_{1r},c_{2r}\in\Z$ and $k,m_r,a_{1r}+a_{2r},b_{1r}+b_{2r},c_{1r}+c_{2r}\in\Z^+$.

{\rm (i)} The tuple $(i+j,i-j,s+t,s-t,u+v,u-v)$ is universal if and only if
$$(a_{1r}+a_{2r},a_{1r}-a_{2r},b_{1r}+b_{2r},b_{1r}-b_{2r},c_{1r}+c_{2r},c_{1r}-c_{2r})$$
is universal for all $r=1,\ldots,k$.

{\rm (ii)} The tuple $(i+j,i-j,s+t,s-t,u+v,u-v)$ is almost universal if and only if the tuple
$$(a_{1r}+a_{2r},a_{1r}-a_{2r},b_{1r}+b_{2r},b_{1r}-b_{2r},c_{1r}+c_{2r},c_{1r}-c_{2r})$$
is almost universal for all $r=1,\ldots,k$.
\end{theorem}

Utilizing the above theorem and some known identities on the theta functions, we establish the following
theorems whose results were conjectured by Sun \cite{Sun2}.

\begin{theorem}\label{Th1.2}
All the tuples
\begin{align*}& (8,2,5,1,3,1),\,(8,2,5,3,4,2),\,(8,2,6,2,3,1),\, (8,2,6,4,4,2),\,(8,2,6,4,6,2),
\\&(8,2,7,1,3,1), (8,4,4,0,3,1),\,(8,4,6,4,4,0),\,(8,6,5,1,3,1),\,(8,6,6,4,4,2),
\\&(8,6,6,4,6,2),\,(8,6,7,1,3,1),\,(9,5,4,2,2,0),\,(9,5,6,2,3,1),\,(9,5,8,4,3,1),\\&
(9,7,6,4,3,1),(9,7,8,4,4,2),\ (9,c,3,1,2,0)\, (c=1,5,7),
\\& (9,c,5,d,3,1)\, (c\in\{1,5,7\},\ d\in\{1,3\}),\,(12,6,6,4,2,0),\, (12,6,6,4,4,2),
\\&(12,6,8,4,6,4),\, (12,6,9,c,3,1)\, (c=1,5,7),\, (12,8,4,2,3,1),\,(12,8,6,4,6,2),
\\&(12,8,9,3,3,1),\,(12,8,12,6,3,1),\,(16,8,6,2,2,0),\,(16,8,8,4,6,4)
\end{align*} are universal. In other words, the sums
\begin{align*}&x(4x+1)+\f{y(3y+1)}2+\f{z(5z+1)}2,\ x(4x+1)+y(2y+1)+\f{z(5z+3)}2,
\\& x(4x+1)+y(3y+1)+\f{z(3z+1)}2,\ x(2x+1)+y(3y+2)+z(4z+1),
\\& x(4x+1)+y(3y+1)+z(3z+2),\ x(4x+1)+\f{y(3y+1)}2+\f{z(7z+1)}2,
\\& 2x^2+y(y+1)+\f{z(3z+1)}2,\ 2x^2+y(y+1)+z(3z+2),
\end{align*}
\begin{align*}
& x(4x+3)+\f{y(3y+1)}2+\f{z(5z+1)}2,\ x(2x+1)+y(3y+2)+z(4z+3),
\\& x(3x+1)+y(3y+2)+z(4z+3),\ x(4x+3)+\f{y(3y+1)}2+\f{z(7z+1)}2,
\\& x^2+\f{y(y+1)}2+\f{z(9z+5)}2,\ x(3x+1)+\f{y(3y+1)}2+\f{z(9z+5)}2,
\\& x(x+1)+\f{y(3y+1)}2+\f{z(9z+5)}2,\  x(3x+2)+\f{y(3y+1)}2+\f{z(9z+7)}2,
\\& x(x+1)+\f{y(y+1)}2+\f{z(9z+7)}2,\ x^2+\f{y(3y+1)}2+\f{z(9z+c)}2\ (c=1,5,7),
\\& \f{x(3x+1)}2+\f{y(5y+d)}2+\f{z(9z+c)}2\ (c\in\{1,5,7\},\ d\in\{1,3\}),
\\&x^2+3y(2y+1)+z(3z+2),\ x(ax+1)+3y(2y+1)+z(3z+2)\,(a=1,2),
\\& 3\f{x(x+1)}2+\f{y(3y+1)}2+\f{z(9z+c)}2\,(c=1,5,7),
\\&\f{x(x+1)}2+\f{y(3y+1)}2+2z(3z+2),\ x(3x+1)+y(3y+2)+2z(3z+2),
\\&2x(3x+2)+\f{y(3y+1)}2+3\f{z(3z+1)}2,\ 3\f{x(x+1)}2+\f{y(3y+1)}2+2z(3z+2),
\\&x^2+2y(y+1)+z(3z+1),\ x(x+1)+2y(y+1)+z(3z+2)
\end{align*}
are all universal.
\end{theorem}

\begin{remark}
It seems that it is almost hopeless to prove Theorem \ref{Th1.2} completely via the theory of ternary quadratic forms. For example, the tuple $(12,8,12,6,3,1)$ is related to the ternary quadratic form
$x^2+9y^2+16z^2$. There are four classes in the genus of $x^2+9y^2+16z^2$, and the three classes not containing $x^2+9y^2+16z^2$ have representatives
$$x^2+y^2+144z^2,\ 2x^2+2y^2+37z^2-2yz+2xz,\ 4x^2+5y^2+10z^2+2yz+4xz+4xy$$
respectively. Also, the tuples $(8,2,7,1,3,1)$ and $(8,6,7,1,3,1)$ are related to the ternary quadratic form $6x^2+14y^2+21z^2$. There are three classes in the genus of $6x^2+14y^2+21z^2$, and the two classes not containing $6x^2+14y^2+21z^2$ have representatives
$$3x^2+14y^2+42z^2\ \ \t{and}\ \ 5x^2+17y^2+21z^2-16yz+4xz+2xy$$
respectively.
\end{remark}

In view of the corollary in \cite[p.\,56]{DSP},
for a particular polynomial $P(x,y,z)=x(a_1x+a_2)+y(b_1y+b_2)+z(c_1z+c_2)\in\Z[x,y,z]$,  it is easy to check that whether $P(x,y,z)$ is almost universal by the theory of ternary quadratic forms.

Using the second part of Theorem \ref{Th1.1}, we deduce the following new result on almost universal sums.

\begin{theorem}\label{Th1.4}  All the sums
\begin{align*}&3x^2+3\f{y(y+1)}2+\f{z(3z+1)}2,\ 3\f{x(x+1)}2+\f{y(3y+1)}2+z(3z+2),
 \\&3x^2+2y(y+1)+\f{z(3z+1)}2,\ 2x(x+1)+\f{y(3y+1)}2+z(3z+2),
 \\&2x^2+3y^2+\f{z(3z+1)}2,\ 2x^2+\f{y(3y+1)}2+z(3z+2),
 \\&x(x+1)+y(3y+1)+\f{z(3z+1)}2,\ x(x+1)+3y(y+1)+\f{z(3z+1)}2,
 \\&2x^2+y(3y+1)+\f{z(3z+1)}2,\ 2x^2+3y(y+1)+\f{z(3z+1)}2
 \end{align*}
are almost universal.
\end{theorem}

We will prove Theorems 1.1-1.3 in Sections 3-5 respectively.

  \section{Universal sums via the theory of quadratic forms} \label{AlmostSec}

  \setcounter{equation}{0}
 \setcounter{conjecture}{0}
 \setcounter{theorem}{0}
 \setcounter{proposition}{0}

  For two integer-valued polynomials $P(x_1,\ldots,x_n)$ and $Q(x_1,\ldots,x_n)$,
  if $$\{P(x_1,\ldots,x_n):\ x_1,\ldots,x_n\in\Z\}=\{Q(x_1,\ldots,x_n):\ x_1,\ldots,x_n\in\Z\}$$
  then we say that $P$ and $Q$ are equivalent and denote this by $P(x_1,\ldots,x_n)\sim Q(x_1,\ldots,x_n)$. For $a_1,b_1,c_1,a_1',b_1',c_1'\in\Z^+$ and $a_2,b_2,c_2,a_2',b_2',c_2'\in\Z$ with $a_1-a_2,b_1-b_2,c_1-c_2,a_1'-a_2',b_1'-b_2',c_1'-c_2'$
 all even, if
 $$\f{x(a_1x+a_2)}2+\f{y(b_1y+b_2)}2+\f{z(c_1z+c_2)}2\sim \f{x(a_1'x+a_2')}2+\f{y(b_1'y+b_2')}2+\f{z(c_1'z+c_2')}2,
$$ then we say that the two tuples $(a_1,a_2,b_1,b_2,c_1,c_2)$ and $(a_1',a_2',b_1',b_2',c_1',c_2')$
are equivalent, and denote this by
  $$(a_1,a_2,b_1,b_2,c_1,c_2)\sim(a_1',a_2',b_1',b_2',c_1',c_2').$$

  As Euler observed,
  \begin{equation}\label{p33}
  p_3(x)+p_3(y) \sim x^2+2p_3(y)
  \end{equation}
  (cf. \cite[p.\,11]{Dickson}). H.-L. Wu and Z.-W. Sun \cite[(1)]{WuSun} proved that
\begin{equation}\label{WuSunequiv}
  p_3(x)+p_5(y)\sim p_5(x)+3p_5(y).
\end{equation}

  \begin{lemma}\label{Lem2.1} {\rm (Sun \cite[Theorem 1.9]{Sun3})} {\rm (i)}
  For any $a\in\Z^+$ and $b\in\N$ with $b\ls a/2$, we have
  \begin{equation}x(ax+b)+y(ay+a-b)\sim ap_3(x)+\f{y(ay+a-2b)}2.
  \end{equation} In particular,
  \begin{equation}\label{2p58}
  2p_5(x)+p_8(y)  \sim 3p_3(x)+p_5(y).
\end{equation}

{\rm (ii)} We have
\begin{align}
x^2+p_3(y) & \sim p_5(x) + 2p_5(y), \label{SunLemma1}\\
p_3(x)+2p_3(y) &\sim p_5(x)+p_8(y),\label{SunLemma2}\\
x^2+4p_3(y) &\sim 4p_5(x)+p_8(y),\label{SunLemma3}
\\p_3(x)+p_3(y)&\sim\f{x(5x+1)}2+\f{y(5y+3)}2.\label{513}
\end{align}
\end{lemma}
\begin{remark} We can also deduce Lemma \ref{Lem2.1}(i) by using the identity
\begin{align}  \label{2.10P}
	f(a, ab^2) f(b, a^2b) =& f(a, b) \psi(ab)
\end{align}
appeared in  \cite[Entry 30 (i)]{Adiga3}.
\end{remark}

  For convenience, for $a,b,c\in\Z^+$ we define
  $$E(a,b,c)=\N\sm\{ax^2+by^2+cz^2:\ x,y,z\in\Z\}.$$
  Dickson \cite[pp.\,112-113]{D39} listed $102$ regular ternary quadratic forms $ax^2+by^2+cz^2$
  for which the structure of $E(a,b,c)$ is known explicitly.

  \begin{theorem}\label{Th2.1} All the tuples
\begin{align*}&(8,2,3,1,2,0),\, (8,2,3,1,3,1),\,(8,2,4,2,2,0),\, (8,2,4,2,4,0),\, (8,2,5,1,4,2),
\\&(8,2,5,3,5,1),\,(8,4,6,4,6,2), \,(8,4,8,2,3,1),\,(8,4,8,4,6,4),\,(8,6,3,1,2,0),
\\&(8,6,4,2,4,0),\, (8,6,4,2,4,2),\,(8,6,5,1,4,2),\,(8,6,8,4,3,1)
\end{align*}
are universal.
\end{theorem}
\Proof.  (a) To prove the universality of $(8,6,4,2,4,2)$
(or the sum $x(4x+3)+y(y+1)/2+z(z+1)/2$), we observe that
\begin{align*} &\ n=x(4x+3)+\f{y(y+1)}2+\f{z(z+1)}2
\\\iff&\ 16n+13=(8x+3)^2+2(2y+1)^2+2(2z+1)^2.
\end{align*}
As $E(1,2,2)=\{4^s(8t+7):\ s,t\in\N\}$ by Dickson \cite[pp.\,112-113]{D39},
for each $n\in\N$ we can write $16n+13=w^2+2u^2+2v^2$ with $w,u,v\in\Z$.
As $2\nmid w$, we have $2(u^2+v^2)\eq13-w^2\eq13-1\eq4\pmod 8$ and hence $2\nmid uv$.
Thus $u=2y+1$ and $v=2z+1$ for some $y,z\in\Z$. Since
$$w^2\eq 13-2(u^2+v^2)\eq 13-2(1+1)=3^2\pmod {16},$$
we have $w\not\eq\pm1\pmod{8}$  and hence we can write $w$ or $-w$ as $8x+3$ with $x\in\Z$.
Thus $16n+13=(8x+3)^2+2(2y+1)^2+2(2z+1)^2$ as desired.

Similarly, $(8,2,4,2,4,2)$ (or $x(2x+1)+y(2y+1)+z(4z+1)$) is universal,
as pointed out in \cite{Sun1}. Since $(8,2,5,3,5,1)\sim (8,2,4,2,4,2)$
by \eqref{513}, the tuple $(8,2,5,3,5,1)$ is also universal.

(b) To prove the universality of $(8,2,3,1,2,0)$
(or the sum $x^2+y(4y+1)+z(3z+1)/2$), we note that
\begin{align*} &\ n=x^2+y(4y+1)+\f{z(3z+1)}2
\\\iff&\ 48n+5=48x^2+3(8y+1)^2+2(6z+1)^2.
\end{align*}
As
$$E(2,3,48)=\bigcup_{k\in\N}\{8k+1,8k+7,16k+6,16k+10,64k+24\}\cup\{9^s(3t+1):\ s,t\in\N\}$$
by Dickson \cite[pp.\,112-113]{D39}, for each $n\in\N$ we can write $48n+5=48x^2+3u^2+2v^2$ with $x,u,v\in\Z$. Clearly, $2\nmid u$ and $2v^2\eq5 -3u^2\eq2\pmod 8$, thus $2\nmid uv$.
As $3u^2\eq 5-2v^2\eq5-2=3\pmod{16}$, we can write $u$ or $-u$ as $8y+1$ with $y\in\Z$.
Since $\gcd(v,6)=1$, we can write $v$ or $-v$ as $6z+1$ with $z\in\Z$. Hence
$48n+5=48x^2+3(8y+1)^2+2(6z+1)^2$ as desired.

Similarly, by using $E(2,3,48)$ we can prove the universality of $(8,6,3,1,2,0)$.

(c) By Dickson \cite[pp.\,112-113]{D39}, we have
\begin{align*} E(2,2,3)&=\{8k+1:\ k\in\N\}\cup\{9^s(9t+6):\ s,t\in\N\},
\\E(1,2,16)&=\bigcup_{k\in\N}\{8k+5,8k+7,16k+10\}\cup\{4^s(16t+14):\ s,t\in\N\},
\\E(1,2,32)&=\bigcup_{k\in\N}\{16k+14,8k+5,2(8k+5),4(8k+5)\}\cup\{4^s(8t+7):\ s,t\in\N\},
\\E(2,3,12)&=\{16k+6:\ k\in\N\}\cup\{9^s(3t+1):\ s,t\in\N\},
\\E(2,3,8)&=\bigcup_{k\in\N}\{8k+1,8k+7,32k+4\}\cup\{9^s(9t+6):\ s,t\in\N\},
\\E(2,5,10)&=\{8k+3:\ k\in\N\}\cup\bigcup_{s,t\in\N}\{25^s(5t+1),25^s(5t+4)\},
\\E(3,3,4)&=\{4k+1:\ k\in\N\}\cup\{9^s(3t+2):\ s,t\in\N\}.
\end{align*}
In view of this, we can easily prove that all the tuples
\begin{align*}&(8,2,3,1,3,1),(8,6,3,1,3,1),(8,2,4,2,2,0),(8,2,4,2,4,0),
(8,6,4,2,4,0),
\\&(8,2,5,1,4,2),(8,6,5,1,4,2),(8,4,6,4,6,2),(8,4,8,2,3,1),(8,4,8,4,6,4),
\\&(8,6,8,4,3,1)
\end{align*}
are universal.
This concludes the proof. \qed

   \begin{lemma}\label{Lem-3} {\rm (Sun \cite[Lemma 3.2]{Sun2015})} Let $w=x^2+3y^2\eq 4\ (\mo\ 8)$ with $x,y\in\Z$.
Then there are odd integers $u$ and $v$ such that $w=u^2+3v^2$.
\end{lemma}

\begin{lemma} \label{Lem-GS} {\rm (F. Ge and Sun \cite[Lemma 2.2]{GeSun1})} Let $w=x^2+3y^2$ with $x,y$ odd and $3\nmid x$.
Then there are integers $u$ and $v$ relatively prime to $6$ such
that $w=u^2+3v^2$.
   \end{lemma}

\begin{theorem}\label{Th2.2} All the tuples
\begin{align*}&(9,1,4,2,3,1),\,(9,3,9,1,3,1),\, (9,5,4,2,4,2),\ (9,5,5,3,5,1),
\\&\,(9,5,8,4,2,0),\,(9,5,9,3,3,1),\,(9,7,4,2,4,2),\, (9,7,5,3,5,1),
\\&  (9,7,8,4,2,0),\ (9,c,3,1,3,1)\, (c=1,5,7),\ (10,6,10,2,6,4)
\end{align*}
are universal.
\end{theorem}
\Proof. In view of \eqref{513} and the universality of $(8,4,8,4,6,4)$ mentioned in Theorem \ref{Th2.1}, we have the universality
of $(10,6,10,2,6,4)$.

To prove the university of $(9,1,3,1,3,1)$ (or $x(3x+1)/2+y(3y+1)/2+z(9z+1)/2$), we note that
\begin{align*}&n=\f{x(3x+1)}2+\f{y(3y+1)}2+\f{z(9z+1)}2
\\\iff&72n+7=3(6x+1)^2+3(6y+1)^2+(18z+1)^2.
\end{align*}
As $E(1,3,3)=\{9^s(3t+2):\ s,t\in\N\}$ by \cite[pp.\,112-113]{D39}, for each $n\in\N$ we can write
$72n+7=w^2+3u^2+3v^2$ with $w,u,v\in\Z$. As $w^2\not\eq 7\pmod4$, $u$ or $v$ is odd.
Without loss of generality, we assume that $v$ is odd. Since $w^2+3u^2\eq 7-3v^2\eq4\pmod8$, by Lemma
\ref{Lem-3} we can write $w^2+3u^2$ as $w_0^2+3u_0^2$ with $w_0$ and $u_0$ odd integers.
Note that $72n+7=w_0^2+3u_0^2+3v^2$. As $3\nmid w_0$, by Lemma \ref{Lem-GS} we can write $w_0^2+3u_0^2$
as $w_1^2+3u_1^2$ with $w_1,u_1\in\Z$ and $\gcd(w_1u_1,6)=1$. As $2\nmid w_1v$ and $3\nmid w_1$,
by Lemma \ref{Lem-GS} we can write $w_1^2+3v^2$ as $w_2^2+3v_1^2$ with $w_2,v_1\in\Z$ and $\gcd(w_2v_1,6)=1$. Write $u_1$ or $-u_1$ as $6x+1$ with $x\in\Z$, and write $v_1$ or $-v_1$
as $6y+1$ with $y\in\Z$. Then $72n+7=w_2^2+3(6x+1)^2+3(6y+1)^2$. As $w_2^2\eq 7-3-3=1\pmod{36}$,
we can write $w_2$ or $-w_2$ as $18z+1$. Thus $72n+1=3(6x+1)^2+3(6y+1)^2+(18z+1)^2$ as desired.

By similar arguments, we can show that both $(9,5,3,1,3,1)$ and $(9,7,3,1,3,1)$ are also universal.

To prove the universality of $(9,1,4,2,3,1)$ (or $x(x+1)/2+y(3y+1)/2+z(9z+1)/2$), we note that
\begin{align*}&n=\f{x(x+1)}2+\f{y(3y+1)}2+\f{z(9z+1)}2
\\\iff&72n+13=9(2x+1)^2+3(6y+1)^2+(18z+1)^2.
\end{align*}
As $E(1,1,3)=\{9^s(9t+6):\ s,t\in\N\}$ by \cite[pp.\,112-113]{D39},
for each $n\in\N$ we can write $72n+13=u^2+v^2+3w^2$ with $u,v,w\in\Z$. As  $3w^2\not\eq13\pmod 4$,
$u$ or $v$ is odd. Without loss of generality, we assume that $2\nmid u$.
As $v^2+3w^2\eq13-u^2\eq4\pmod8$, by Lemma \ref{Lem-3} we may write $v^2+3w^2=v_1^2+3w_1^2$
with $v_1$ and $w_1$ odd integers. Thus $72n+13=u^2+v_1^2+3w_1^2$ and $2\nmid uv_1w_1$.
As $u^2+v_1^2\eq13\eq1\pmod3$, $u$ or $v_1$ is divisible by $3$. Without loss of generality, we assume that $u=3(2x+1)$ with $x\in\Z$. Since $2\nmid v_1w_1$ and $3\nmid v_1$, by Lemma \ref{Lem-GS} we can write
$v_1^2+3w_1^2$ as $v_2^2+3w_2^2$ with $v_2,w_2\in\Z$ and $\gcd(v_2w_2,6)=1$.
We may write $w_2$ or $-w_2$ as $6y+1$ with $y\in\Z$. Thus
$72n+13=9(2x+1)^2+3(6y+1)^2+v_2^2$. As $v_2^2\eq 13-9-3=1\pmod{36}$, we may write $v_2$ or $-v_2$
as $18z+1$ with $z\in\Z$. So $72n+13=9(2x+1)^2+3(6y+1)^2+(18z+1)^2$ as desired.

In view of \eqref{WuSunequiv}, from the universality of
$(9,1,4,2,3,1)$, we deduce that the tuple $(9,3,9,1,3,1)$ is universal.
The tuple $(9,5,4,2,3,1)$ is known to be universal (cf. \cite[Theorem 1.2(ii)]{Sun2}),
so $(9,5,9,3,3,1)$ is also universal with the aid of \eqref{WuSunequiv}.

To prove the universality of $(9,5,4,2,4,2)$ (or $x(x+1)/2+y(y+1)/2+z(9z+5)/2$), we note that
\begin{align*}&n=\f{x(x+1)}2+\f{y(y+1)}2+\f{z(9z+5)}2
\\\iff& 72n+43=9(2x+1)^2+9(2y+1)^2+(18z+5)^2.
\end{align*}
As $E(1,1,1)=\{4^s(8t+7):\ s,t\in\N\}$ (cf. \cite[pp.\,112-113]{D39}, for each $n\in\N$ we can write
$72n+43=u^2+v^2+w^2$ with $u,v,w\in\Z$. As $u^2+v^2+w^2\eq43\eq3\pmod8$, we have $2\nmid uvw$.
Since $u^2+v^2+w^2\eq1\pmod3$, without loss of generality we may assume that both $u$ and $v$
are divisible by $3$. Thus  $u=3(2x+1)$ and $v=3(2y+1)$ for some $x,y\in\Z$.
As $w^2\eq43-9-9=5^2\pmod{36}$, we may write $w$ or $-w$ as $18z+5$ with $z\in\Z$.
So $72n+43=9(2x+1)^2+9(2y+1)^2+(18z+5)^2$ as desired.

Similarly, we can prove that $(9,7,4,2,4,2)$ is universal.

For $c=5,7$ we have
$$ (9,c,5,3,5,1)\sim(9,c,4,2,4,2)\sim(9,c,8,4,2,0)$$
by \eqref{513} and \eqref{p33}. Since $(9,5,4,2,4,2)$ and $(9,7,4,2,4,2)$ are universal, the tuples
$$(9,5,5,3,5,1),\ (9,5,8,4,2,0),\ (9,7,5,3,5,1),\ (9,7,8,4,2,0)$$
are also universal.

In view of the above, we have completed the proof of Theorem \ref{Th2.2}. \qed

The following known lemma can be found in \cite[Lemma 2.1]{Sun2015}.

  \begin{lemma}\label{Lem-m}  Let $w=x^2+my^2$ be a positive integer
with $m\in\{2,5,8\}$ and $x,y\in\Z$. Then we can write $w$ in the form $u^2+mv^2$ with $u,v\in\Z$
such that $u$ or $v$ is not divisible by $3$.
  \end{lemma}

  The following lemma follows from \cite[Lemma 2.2(i)]{S16} (or the comments around \cite[(1.4)]{Sun2015}).

  \begin{lemma}\label{Lem-9} If $x,y,z\in\Z$ are not all divisible by $3$, then
  $9(x^2+y^2+z^2)=u^2+v^2+w^2$ for some $u,v,w\in\Z$ not all divisible by $3$.
  \end{lemma}

 \begin{theorem}\label{Th2.3} All the tuples
\begin{align*}&(12,4,3,1,2,0),\,(12,4,6,4,6,2), \,(12,8,3,1,3,1),
\\&(12,6,8,2,3,1),\, (12,6,8,6,3,1),\,(12,8,3,1,2,0),\,(12,8,5,3,5,1),
\\&(12,8,4,2,4,2),\,(12,8,8,4,2,0),\,(12,8,8,4,3,1),\,(12,8,12,4,2,0)
\end{align*}
are universal.
\end{theorem}
\Proof. To prove the universality of $(12,4,3,1,2,0)$ (or $x^2+2y(3y+1)+z(3z+1)/2$), we observe that
\begin{align*}&n=x^2+2y(3y+1)+\f{z(3z+1)}2
\\\iff&24n+5=24x^2+4(6y+1)^2+(6z+1)^2.
\end{align*}
By Dickson \cite[pp.\,112-113]{D39},
$$E(1,4,24)=\bigcup_{k\in\N}\{4k+2,4k+3\}\cup\{9^s(9t+3):\ s,t\in\N\}.$$
Thus, for each $n\in\N$ we can write $24n+5=24x^2+4u^2+v^2$ with $x,u,v\in\Z$.
Clearly, $2\nmid v$. As $4u^2\eq5-v^2\eq4\pmod8$, $u$ is also odd. Since
$u^2+v^2\eq 4u^2+v^2\eq5\eq2\pmod3$, we have $3\nmid uv$. As $uv$ is relatively prime to $6$,
we can write $u$ or $-u$ as $6y+1$, and write $v$ or $-v$ as $6z+1$, where $y$ and $z$ are integers.
Now, $n=24x^2+4(6y+1)^2+(6z+1)^2$ as desired.

To prove the universality of $(12,6,8,2,3,1)$ (or $3x(x+1)/2+y(3y+1)/2+z(4z+1)$), we note that
\begin{align*}&n=3\f{x(x+1)}2+\f{y(3y+1)}2+z(4z+1)
\\\iff&48n+23=18(2x+1)^2+2(6y+1)^2+3(8z+1)^2.
\end{align*}
By Dickson \cite[pp.\,112-113]{D39},
$$E(2,3,18)=\bigcup_{k\in\N}\{3k+1,8k+1\}\cup\{9^s(9t+3):\ s,t\in\N\}.$$
Thus, for each $n\in\N$ we can write $48n+23=18u^2+2v^2+3w^2$ with $u,v,w\in\Z$.
Clearly, $2\nmid w$. Since $2(u^2+v^2)\eq 23-3w^2\eq4\pmod{8}$, we have $2\nmid uv$.
Note also that $3\nmid v$. So we may write $u=2x+1$ with $x\in\Z$, and $v$ or $-v$ as $6y+1$
with $y\in\Z$. As $3w^2\eq 23-2u^2-2v^2\eq 23-2-2\eq3\pmod{16}$, we may write $w$ or $-w$ as
$8z+1$ with $z\in\Z$. Thus $48n+23=18(2x+1)^2+2(6y+1)^2+3(8z+1)^2$ as desired.
Similarly, $(12,6,8,6,3,1)$ is universal.

By similar arguments, we can prove that $(12,8,3,1,2,0)$ is universal by using
$$E(1,16,24)=\bigcup_{k\in\N}\{4k+2,4k+3,8k+5,64k+8\}\cup\{9^s(9t+3):\ s,t\in\N\}$$
(cf. \cite[pp.\,112-113]{D39}).
We can also prove that $(12,8,8,4,3,1)$ is universal by using
$$E(1,6,16)=\bigcup_{k\in\N}\{8k+3,8k+5,16k+2,16k+14,64k+8\}\cup\{9^s(9t+3):\ s,t\in\N\}$$
(cf. \cite[pp.\,112-113]{D39}), and that the three tuples
$$(12,8,5,3,5,1),\ (12,8,4,2,4,2) \ \t{and}\ (12,8,8,4,2,0)$$
(which are equivalent in view of \eqref{p33} and \eqref{513})
are universal by using
$$E(3,8,12)=\bigcup_{k\in\N}\{4k+1,4k+2\}\cup\{9^s(3t+1):\ s,t\in\Z\}$$
(cf. \cite[pp.\,112-113]{D39}).
Also, we can prove that $(12,8,12,4,2,0)$ is universal by using
$$E(1,4,6)=\{16k+2:\ k\in\N\}\cup\{9^s(9t+3):\ s,t\in\N\}$$
(cf. \cite[pp.\,112-113]{D39}).

To prove the universality of $(12,4,6,4,6,2)$ (or $2x(3x+1)+y(3y+1)+z(3z+2)$), we note that
\begin{align*}&n=2x(3x+1)+y(3y+1)+z(3z+2)
\\\iff&12n+7=2(6x+1)^2+(6y+1)^2+4(3z+1)^2.
\end{align*}
As $E(1,2,4)=\{4^s(16t+14):\ s,t\in\N\}$ (cf. \cite[pp.\,112-113]{D39}),
for each $n\in\N$ we can write $12n+7$ as $u^2+2v^2+4w^2=u^2+2(v^2+2w^2)$ with $u,v,w\in\Z$.
By Lemma \ref{Lem-m}, without loss of generality, we may assume that $u$ or $v$ is not divisible by $3$, and also $v$ or $w$ is not divisible by $3$. If $3\mid v$, then $3\nmid uw$ and also
$$2\eq u^2+w^2\eq u^2+4w^2\eq7\eq1\pmod 3.$$
Therefore we must have $3\nmid v$. As $u^2+w^2\eq u^2+4w^2\eq 7-2v^2\eq2\pmod3$, we also have $3\nmid uw$.
As $u$ is relatively prime to $6$, we may write $u$ or $-u$ as $6y+1$ with $y\in\Z$.
Note that $2v^2\eq7-u^2\eq2\pmod4$. So we may write $v$ or $-v$ as $6x+1$ with $x\in\Z$.
Clearly, we may write $w$ or $-w$ as $3z+1$ with $z\in\Z$.  Therefore,
$12n+7=2(6x+1)^2+(6y+1)^2+4(3z+1)^2$ as desired.

To prove the universality of $(12,8,3,1,3,1)$ (or $2x(3x+2)+y(3y+1)/2+z(3z+1)/2$), we note that
\begin{align*}&n=2x(3x+2)+\f{y(3y+1)}2+\f{z(3z+1)}2
\\\iff&24n+18=16(3x+1)^2+(6y+1)^2+(6z+1)^2.
\end{align*}
As $E(1,1,1)=\{4^s(8t+7):\ s,t\in\N\}$ (cf. \cite[pp.\,112-113]{D39}),
for each $n\in\N$ we can write $24n+18$ as a sum of three squares.
In view of Lemma \ref{Lem-9}, we may write $24n+18=u^2+v^2+w^2$ with $u,v,w\in\Z$ not all divisible by $3$. As $u^2+v^2+w^2\eq0\pmod3$, we must have $3\nmid uvw$. Without loss of generality, we may assume that $w$ is even. As $u^2+v^2\eq 18\eq2\pmod4$, we have $2\nmid uv$. Note that $w^2\eq 18-u^2-v^2\eq0\pmod 8$
and hence $w\eq0\pmod4$. We may write $w$ or $-w$ as $4(3x+1)$ with $x\in\Z$. Also, we may write $u$ or $-u$ as $6y+1$ with $y\in\Z$, and write $v$ or $-v$ as $6z+1$ with $z\in\Z$. Thus
$24n+18=16(3x+1)^2+(6y+1)^2+(6z+1)^2$ as desired.
This ends the proof. \qed

\begin{theorem}\label{Th2.4} All the tuples
\begin{align*}&(16,8,6,4,4,0),\,(16,8,6,4,4,2),\,(18,6,6,4,6,2),
\\&(20,10,5,1,4,2),\,(20,10,5,3,4,2),\,(24,12,6,2,2,0)
\end{align*} are universal.
\end{theorem}
\Proof. In view of \eqref{p33}, we have $(16,8,6,4,4,0)\sim(8,4,8,4,6,4)$. As $(8,4,8,4,6,4)$
is universal by Theorem \ref{Th2.1}, the tuple $(16,8,6,4,4,0)$ is also universal.

To prove the universality of $(16,8,6,4,4,2)$ (or $x(x+1)/2+2y(y+1)+z(3z+2)$), we note that
\begin{align*}&n=\f{x(x+1)}2+2y(y+1)+z(3z+2)
\\\iff&24n+23=3(2x+1)^2+12(2y+1)^2+8(3z+1)^2.
\end{align*}
As
$$E(3,8,12)=\bigcup_{k\in\N}\{4k+1,4k+2\}\cup\{9^s(3t+1):\ s,t\in\N\}$$
by \cite[pp.\,112-113]{D39}, for each $n\in\N$ we can write $24n+23$
as $3u^2+12v^2+8w^2$ with $u,v,w\in\Z$. Clearly $u$ is odd.
As $12v^2\eq23-3u^2\eq 20\not\eq0\pmod{8}$, $v$ is odd. So $u=2x+1$ and $v=2y+1$ for some $x,y\in\Z$.
Since $8w^2\eq 23\pmod3$, we can write $w$ or $-w$ as $3z+1$ with $z\in\Z$.
Therefore $24n+23=3(2x+1)^2+12(2y+1)^2+8(3z+1)^2$ as desired.

By \eqref{WuSunequiv}, we have $(18,6,6,4,6,2)\sim(8,4,6,4,6,2)$. Since $(8,4,6,4,6,2)$
is universal by Theorem \ref{Th2.1}, the tuple $(18,6,6,4,6,2)$ is also universal.

 To prove the universality of $(20,10,5,1,4,2)$ (or the sum $x(x+1)/2+5y(y+1)/2+z(5z+1)/2$), we note that
\begin{align*}&\ n=\f{x(x+1)}2+5\f{y(y+1)}2+\f{z(5z+1)}2
\\\iff&\ 40n+31=5(2x+1)^2+25(2y+1)^2+(10z+1)^2.
\end{align*}
By Dickson \cite[pp.\,112-113]{D39},
$$E(1,5,25)=\bigcup_{k\in\N}\{5k+2,5k+3,25k+10,25k+15\}\cup\{4^s(8t+3):\ s,t\in\N\}.$$
Thus, for any $n\in\N$ we can write $40n+31$ as $w^2+5u^2+25v^2$ with $w,u,v\in\Z$.
If $w$ is even, then
$$u^2+v^2\eq 5u^2+25v^2\eq 31-w^2\eq3\pmod 4$$
which is impossible. Note also that $w^2\eq31\eq1\pmod5$ and hence $w\eq\pm1\pmod5$.
Thus, $w$ or $-w$ can be written as $10z+1$ with $z\in\Z$. As $u^2+v^2\eq 31-w^2\eq2\pmod4$,
we have $2\nmid uv$. Thus $u=2x+1$ and $v=2y+1$ for some $x,y\in\Z$. Therefore
$$40n+31=w^2+5u^2+25v^2=5(2x+1)^2+25(2y+1)^2+(10z+1)^2$$
as desired.

Similarly, we can prove that the tuple $(20,10,5,3,4,2)$ is universal.

In view of \eqref{2p58}, we have $(24,12,6,2,2,0)\sim(12,8,12,4,2,0)$.
As $(12,8,12,4,2,0)$ is universal by Theorem \ref{Th2.3}, the tuple $(24,12,6,2,2,0)$
is also universal.

In view of the above, we have completed the proof of Theorem \ref{Th2.4}. \qed

For a  tuple $(a_1,a_2,b_1,b_2,c_1,c_2)$ with related ternary quadratic form having class number $2$, occasionally we may employ the following lemma to prove its universality.

\begin{lemma}\label{Lem-genus} {\rm (\cite[Theorem 1.3]{C})} Let $f$ be an integral quadratic form with nonzero discriminant. If an integer $m$ is represented by $f$ over the field of real numbers
as well as the ring $\Z_p$ of $p$-adic integers for each prime $p$, then $m$ is represented over $\Z$ by some form $f^*$ in the same genus as $f$.
\end{lemma}

For example, the tuples $(20,10,5,1,3,1)$ and $(20,10,5,3,3,1)$
were proved to be universal in \cite[Theorem 1.1(iv)]{Sun3} with the aid of Lemma \ref{Lem-genus}.
There are two classes in the genus of the related quadratic form $3x^2+3y^2+5z^2$, the one not containing $3x^2+3y^2+5z^2$ has the representative
$3x^2+2y^2+8z^2-2yz$.

\section{An approach via Ramanujan's theta functions}\label{UnivSec}

Now we turn to a new approach to universal sums via Ramanujan's theta fundtions.

\medskip
\noindent{\tt Proof of Theorem \ref{Th1.1}}. Let $|q|<1$. For $a,b\in\Z$ with $a+b>0$, clearly
$|q^aq^b|<1$ and
$$f(q^a,q^b)=\sum_{x=-\infty}^\infty q^{ax(x+1)/2+bx(x-1)/2}=\sum_{x=-\infty}^\infty q^{x((a+b)x+a-b)/2}.$$

For any $n\in\N$, let $R(n)$ denote the number of triples $(x,y,z)\in\Z^3$ satisfying the equation
$$\f{x((i+j)x+i-j)}2+\f{y((s+t)y+s-t)}2+\f{z((u+v)z+u-v)}2=n,$$
and let $R_r(n)$ (with $1\ls r\ls k$) denote the number of ways to write $n$ as
$$ \f{x((a_{1r}+a_{2r})x+a_{1r}-a_{2r})}2+\f{y((b_{1r}+b_{2r})y+b_{1r}-b_{2r})}2
+\f{z((c_{1r}+c_{2r})z+c_{1r}-c_{2r})}2$$
with $(x,y,z)\in\Z^3$.

For any fixed $n\in\N$ and $r\in\{1,\ldots,k\}$, by comparing coefficients of $q^{kn+r-1}$
on both sides of \eqref{=}, we obtain
\begin{equation}\label{kn+r}R(kn+r-1)=m_rR_r(n).
\end{equation}
Thus, for any $n_0\in\N$, we have
\begin{equation}\label{n0}\begin{aligned} & R_r(n)>0\ \t{for all}\ r=1,\ldots,k\ \t{and}\ n=n_0,n_0+1,\ldots
\\ &\ \ \iff  R(n)>0\ \t{for every integer}\ n\gs kn_0.
\end{aligned}
\end{equation}
When $n_0=0$, this yields part (i) of Theorem \ref{Th1.1}.
As \eqref{n0} holds for any $n_0\in\N$, we also have part (ii) of Theorem \ref{Th1.1}.
This concludes the proof. \qed

\begin{remark} In view of \eqref{kn+r}, $|\{n\in\N:\ R(n)=0\}|=1$, if and only if for certain $1\ls r_0\ls k$ we have
$|\{n\in\N:\ R_{r_0}(n)=0\}|=1$ and $R_r(n)>0$ for all $r\in\{1,\ldots,k\}\sm\{r_0\}$ and $n\in\N$.
\end{remark}

In $2023$, Ju \cite{Ju} classified  all almost universal sums of triangular numbers with one exception.
Furthermore, he provided an effective criterion for almost universality
with one exception of an arbitrary sum of triangular numbers, which might be considered as a
natural generalization of the 15-theorem of Conway, Miller and Schneeberger (cf. \cite{CGSS}).

To apply Theorem \ref{Th1.1}, we need to make certain preparation.

Let $a$ and $b$ be complex numbers with $|ab|<1$. The Ramanujan theta function $f(a,b)$ satisfies the following
basic properties \cite{Adiga3}:
\begin{align}\label{2.3}f(a,b)& = f(b,a)\\
\intertext{and}
\label{2.4}f(1,a)& = 2f(a,a^{3}).
\end{align}
The functions $\varphi(q)$ given by \eqref{2.7}, $\psi(q)$ given by \eqref{2.8},
\begin{equation}\label{X(q)}
 X(q):= f(q,q^2) =  \sum_{x=-\infty}^{\infty} {q^{x(3x+1)/2}}
 \end{equation}
 and
\begin{equation}\label{Y(q)} Y(q):= f(q,q^5) =  \sum_{x=-\infty}^{\infty} {q^{x(3x+2)}}
\end{equation}
are related to the squares, the triangular numbers, the generalized  pentagonal numbers and the generalized octagonal  numbers, respectively.
In general, the  theta function $f(q, q^{m-3})$ with $m\gs3$
is related to the generalized $m$-gonal numbers  $p_m(x)$ $(x\in\Z)$.

Let $|ab|<1$. The following general formula is due to  Ramanujan (cf. \cite[ Entry 31]{Adiga3}:
\begin{equation}\label{RamIden}
  f(U_1,V_1)=\sum_{r=0}^{n-1}U_r f \left( \frac{U_{n+r}}{U_r},\frac{V_{n-r}}{U_r}  \right),
\end{equation}
where
$$U_k=a^{k(k+1)/2}b^{k(k-1)/2}\ \ \t{and}\ \ V_k=a^{k(k-1)/2}b^{k(k+1)/2}.$$
This result implies the following lemma.

\begin{lemma} \label{Lem-3.1} For $|q|<1$, we have
\begin{align}\label{appli12}
  \varphi(q)&=\varphi\left(q^4\right)+2q\psi\left(q^8\right),\\
  \label{AAthm16.3}  \psi(q)&=f\left(q^6,q^{10}\right)+qf\left(q^2,q^{14}\right),\\
  \label{Athm20u1}  Y(q)&=X\left(q^8\right)+qY\left(q^4\right),\\
  \label{AAthm216aa} X(q)&= f\left(q^{12},q^{15}\right)+qf\left(q^{6},q^{21}\right)+q^2f\left(q^{3},q^{24}\right).
\end{align}
\end{lemma}

Let $|q|<1$. In \cite[(50)]{Adiga2}, setting $k=3$, $l=2$ and $\varepsilon_1=\varepsilon_2=1$,  we obtain
\begin{equation}\label{aAthm13u1}
  X(q)X(q^2)=\varphi(q^9)X(q^3)+qX(q^3)Y(q^3)+2q^2\psi(q^9)Y(q^3).
\end{equation}
 In \cite[(15)]{Adiga2}, setting $k=3$,  $l=2$, $g=5$, $h=1$, $u=2$, $v=1$ and $\varepsilon_1=\varepsilon_2=1$,  we obtain
\begin{equation}\label{aAthm13u6}
  X(q)Y(q)=X(q^3)X(q^6)+2q\psi(q^9)X(q^6)+2q^2\psi(q^{18})X(q^3).
\end{equation}

Lemma \ref{Lem-3.1} and the identities \ref{aAthm13u1} and \ref{aAthm13u6} will be used in the next two sections.

\section{Proof of Theorem \ref{Th1.2}}\label{Sec4}

\medskip
\noindent{\tt Proof of Theorem \ref{Th1.2}}.
Let $|q|<1$. Multiplying both sides of the identity \eqref{AAthm16.3} by $f(q^6,q^4)f(q^4,q^2)$, we get
  \begin{equation}\label{Athm22u2}
     f(q^6,q^4)f(q^4,q^2)f(q^3,q)=f(q^{10},q^6)f(q^6,q^4)f(q^4,q^2)+qf(q^{14},q^2)f(q^6,q^4)f(q^4,q^2).
  \end{equation}
Recall that the tuple  $(10,2,6,2,4,2)$ (or the sum  $x(2x+1)+y(3y+1)+z(5z+1)$) is universal as
conjectured by Sun \cite{Sun1} and confirmed by
  Ju and Oh \cite{JuOh}. Thus, by applying  Theorem \ref{Th1.1}(i) to the identity
  \eqref{Athm22u2}, we find that
the tuples $(8,2,5,1,3,1)$ and $(8,6,5,1,3,1)$ are universal.

Using \eqref{AAthm16.3}, we find that
  \begin{equation}\label{Athm22u21}
 \psi(q)\psi(q^2)f(q^2,q^8)=f\left(q^6,q^{10}\right)\psi(q^2)f(q^2,q^8)+qf\left(q^2,q^{14}\right)\psi(q^2)f(q^2,q^8).
\end{equation}
Sun \cite{Sun2015} determined the universality of the tuple $(10,6,8,4,4,2)$ (or the sum $p_3+2p_3+2p_7$). Thus, using the first term of the identity
  \eqref{Athm22u21} with the help of Theorem \ref{Th1.1}(i), we see that the tuple $(8,2,5,3,4,2)$ is universal.
  From \eqref{AAthm16.3}, we obtain
  \begin{equation}\label{Athm22u3}
     \psi(q)X\left(q^2\right)f\left(q^6,q^{8}\right)=X\left(q^2\right)f\left(q^6,q^{10}\right)f\left(q^6,q^{8}\right)+qX\left(q^2\right)f\left(q^2,q^{14}\right)f\left(q^6,q^{8}\right).
  \end{equation}
   Ju and Oh \cite{JuOh} proved that the sum $x(2x+1)+y(3y+1)+z(7z+1)$ (or the tuple $(14,2,6,2,4,2)$) is universal as conjectured by Sun \cite{Sun1}.
  Thus the identity  \eqref{Athm22u3} with the help of Theorem \ref{Th1.1}(i), yields the universality of the tuples $(8,2,7,1,3,1)$ and $(8,6,7,1,3,1)$.

  By \eqref{AAthm16.3}, we have
  \begin{equation}\label{univ1}
   \psi(q)X(q^2)X(q^4)=X(q^2)X(q^4)f\left(q^6,q^{10}\right)+qX(q^2)X(q^4)f\left(q^2,q^{14}\right).
  \end{equation}
  Sun \cite{Sun} proved that the tuple $(8,4,4,2,4,0)$ is  universal. As $(8,4,4,2,4,0) \sim (12,4,6,2,4,2)$ by \eqref{SunLemma1}, the tuple  $(12,4,6,2,4,2)$ is universal.
  Combining this with the first term of \eqref{univ1} and Theorem \ref{Th1.1}(i), we obtain the universality of the tuple $(8,2,6,2,3,1)$.

  Using \eqref{AAthm16.3}, we find that
 \begin{equation}\label{Athm16.8au}
   \psi(q)\psi\left(q^2\right)Y\left(q^2\right)=\psi\left(q^2\right)Y\left(q^2\right)f\left(q^{6},q^{10}\right)+q\psi\left(q^2\right)Y\left(q^2\right)f\left(q^2,q^{14}\right).
 \end{equation}
As conjectured by Sun \cite[Conjecture 1.13]{Sun2015} and confirmed by
Ju, Oh and Seo \cite{JuOhSeo}, the tuple $(12,8,8,4,4,2)$  is universal.
Combining this with the identity  \eqref{Athm16.8au} and Theorem \ref{Th1.1}(i), we get the universality of the tuples $(8,2,6,4,4,2)$ and $(8,6,6,4,4,2)$.

  In view of \eqref{AAthm16.3},
\begin{equation}\label{Athm16.10au}
  \psi(q)X\left(q^4\right)Y\left(q^2\right)=X\left(q^4\right)Y\left(q^2\right)f\left(q^{6},q^{10}\right)+qX\left(q^4\right)Y\left(q^2\right)f\left(q^2,q^{14}\right).
 \end{equation}
  It is known (cf. \cite{Sun2}) that the tuple $(24, 12, 6, 2, 4, 2)$ is universal.
   By \eqref{2p58}, $(24, 12, 6, 2, 4, 2) \sim (12,8,12,4,4,2)$. Thus the tuple $(12,8,12,4,4,2)$ is also universal.
  Combining this with \eqref{Athm16.10au} and Theorem \ref{Th1.1}(i), we obtain the universality of the tuples $(8,2,6,4,6,2)$ and $(8,6,6,4,6,2)$.

Using \eqref{appli12}, we find that
\begin{equation}\label{Athm23u6} \begin{aligned}
 \varphi\left(q\right) \varphi\left(q^{2}\right) X\left(q^{4}\right)=&\ \varphi\left(q^4\right) \varphi\left(q^{8}\right) X\left(q^{4}\right)+2q \varphi\left(q^8\right) \psi\left(q^{8}\right) X\left(q^{4}\right) \\
 &\ +2q^2\varphi\left(q^4\right) \psi\left(q^{16}\right) X\left(q^{4}\right)+4q^3\psi\left(q^8\right) \psi\left(q^{16}\right) X\left(q^{4}\right).
 \end{aligned}
 \end{equation}
By Sun \cite[Theorem 1.2(i)]{Sun2}, the tuple $(12,4,4,0,2,0)$ is universal. Thus, by
 applying Theorem \ref{Th1.1}(i) on the identity \eqref{Athm23u6}, we see that the second term and the fourth term yield the universality of the tuples
  $(8,4,4,0,3,1)$ and $(16,8,8,4,3,1)$.

From  \eqref{Athm20u1}, we have
\begin{equation}\label{univ3}
   \varphi\left(q^2\right) \psi\left(q^{2}\right) Y(q)= \varphi\left(q^2\right) \psi\left(q^{2}\right)X\left(q^8\right)+q \varphi\left(q^2\right) \psi\left(q^{2}\right)Y\left(q^4\right).
\end{equation}
Since the tuples $(12,4,4,2,2,0)$ and $(12,8,4,2,2,0)$ are universal (cf. \cite{Sun2}),
by applying Theorem \ref{Th1.1}(i) to the identity  \eqref{univ3}, we obtain the universality of the
  tuple $(8,4,6,4,4,0)$.

 From \eqref{AAthm216aa}, we obtain
    \begin{align}\label{Athm19u3}
  \varphi \left(q^3\right) \psi \left(q^3\right)  X(q)=&\ \varphi \left(q^3\right) \psi \left(q^3\right)f\left(q^{12},q^{15}\right)+q\varphi \left(q^3\right) \psi \left(q^3\right)f\left(q^{6},q^{21}\right) \nonumber \\&\ +q^2\varphi \left(q^3\right) \psi \left(q^3\right)f\left(q^{3},q^{24}\right).
  \end{align}
   As the sum $3x^2+y(3y+1)+z(3z+2)$  (or the tuple $(6,4,6,2,6,0)$) is universal by Sun  \cite{Sun1}, with the aid of \eqref{2p58} the tuple $(12,6,6,0,3,1)$ is also universal. Thus, using the second term of \eqref{Athm19u3} with the help of Theorem \ref{Th1.1}(i), we get the universality of the tuple $(9,5,4,2,2,0)$.

 In a similar way, using \eqref{AAthm216aa} and Theorem \ref{Th1.1}(i), we determine the universality of the tuples listed in the last column of the following table:
  \begin{center}
\begin{tabular}{|c|c|c|c|c|}
\hline
 The product & The tuple used & found in&term/s used & universal tuple/s  \\ \hline
 &&&&$(9,1,5,1,3,1)$ \\
 $ X(q)X\left(q^3\right)f\left(q^9,q^6\right)  $& $(15, 3, 9, 3, 3, 1)$ & \cite{Sun2}&all terms&$(9,5,5,1,3,1)$ \\
 &&&&$(9,7,5,1,3,1)$ \\ \hline
 &&&&$(9,1,5,3,3,1)$ \\
 $ X(q)X\left(q^3\right)f\left(q^{12},q^3\right)  $& $(15, 9, 9, 3, 3, 1)$ & \cite{Sun2}&all terms&$(9,5,5,3,3,1)$ \\
 &&&&$(9,7,5,3,3,1)$ \\ \hline
$X(q) X\left(q^3\right) X \left(q^6\right) $& $(18,6,9,3,3,1)$  &\cite{Guy}  & second term  & $(9,5,6,2,3,1)$\\ \hline
 $ \psi \left(q^6\right) X(q) X\left(q^3\right)  $& $(24,12,9,3,3,1)$ &  \cite{Sun2} & second term  & $(9,5,8,4,3,1)$ \\ \hline
 $ X(q) X\left(q^3\right) Y\left(q^3\right) $& $(18,12,9,3,3,1)$  & \cite{Sun2} & third term  & $(9,7,6,4,3,1)$ \\ \hline
 $\psi \left(q^3\right) \psi \left(q^6\right)  X(q)$& $(24,12,12,6,3,1)$& \cite{Sun2} & third term  & $(9,7,8,4,4,2)$ \\ \hline
 &&&&$(9,1,3,1,2,0)$ \\
 $ \varphi \left(q^3\right) X(q) X\left(q^3\right)  $& $(9,3,6,0,3,1)$ & \cite{Sun2}&all terms&$(9,5,3,1,2,0)$ \\
 &&&&$(9,7,3,1,2,0)$ \\ \hline
 &&&&$(12,6,9,1,3,1)$ \\
 $ \psi \left(q^9\right) X(q) X\left(q^3\right)  $& $(36,18,9,3,3,1)$ &  \cite{JuOhSeo}&all terms&$(12,6,9,5,3,1)$ \\
 &&&&$(12,6,9,7,3,1)$ \\ \hline
  \end{tabular}
\end{center}
\medskip

The universality of $(36,18,9,3,3,1)$ in the above table follows from
 the universality of the tuple $(36,18,4,2,3,1)$ (or the sum $p_3+9p_3+p_5$)
 (which was conjectured by \cite[Conjecture 1.13]{Sun2015} and confirmed by Ju, Oh and Seo \cite{JuOhSeo}, since $(36,18,4,2,3,1) \sim (36,18,9,3,3,1)$ by \eqref{WuSunequiv}.

Using \eqref{aAthm13u1}, we have
\begin{equation}\label{aAthm13u2}
 \varphi(q^3) X(q)X(q^2)=\varphi(q^3)\varphi(q^9)X(q^3)+q\varphi(q^3)X(q^3)Y(q^3)+2q^2\varphi(q^3)\psi(q^9)Y(q^3).
\end{equation}
Sun \cite[Section 3]{Sun2} showed  the universality of  the tuple $(6, 2, 6, 0, 3, 1)$ over $\mathbb{Z}$. Combining this with the third term of \eqref{aAthm13u2} and Theorem \ref{Th1.1}(i), we get the universality of the tuple $(12,6,6,4,2,0)$.

Again, using \eqref{aAthm13u1}, we find that
\begin{equation}\label{aAthm13u3}
 \psi(q^3) X(q)X(q^2)=\varphi(q^9)\psi(q^3)X(q^3)+q\psi(q^3)X(q^3)Y(q^3)+2q^2\psi(q^3)\psi(q^9)Y(q^3).
\end{equation}
By \cite{GPS}, $x^2+p_3(y)+3p_3(z)$ is universal.
So the tuple $(12, 6, 6, 2, 3, 1)$ is universal over $\mathbb{Z}$ with the aid of \eqref{SunLemma1}.
Combining this with the third term of the identity  \eqref{aAthm13u3} and Theorem \ref{Th1.1}(i), we obtain the universality of the tuple $(12,6,6,4,4,2)$.

 From \eqref{aAthm13u1}, we obtain
\begin{equation}\label{aAthm13u4}
 \psi(q^6) X(q)X(q^2)=\varphi(q^9)\psi(q^6)X(q^3)+q\psi(q^6)X(q^3)Y(q^3)+2q^2\psi(q^6)\psi(q^9)Y(q^3).
\end{equation}
As $x^2+p_3(y)+6p_3(z)$ is universal by \cite{GPS}, the tuple $(24, 12, 6, 2, 3, 1)$ is universal
in view of \eqref{SunLemma1}.
This, together with the third term of the identity  \eqref{aAthm13u4} and Theorem \ref{Th1.1}(i), yields the universality of the tuple $(12,6,8,4,6,4)$.

Using \eqref{Athm20u1}, we find that
 \begin{equation}\label{Athm23u8}
 \psi\left(q^{2}\right) X\left(q^{2}\right) Y(q)= \psi\left(q^{2}\right) X\left(q^{2}\right) X\left(q^8\right)+q  \psi\left(q^{2}\right) X\left(q^{2}\right) Y\left(q^4\right).
\end{equation}
By \cite[Theorem 1.11]{Sun2015}, the tuple $(12, 6, 8, 4, 3, 1)$ is universal over $\mathbb{Z}$.
In view of \eqref{2p58}, this is equivalent to the university of
  $(8,4,6,4,6,2)$. Combining this with the second term of \eqref{Athm23u8} and Theorem \ref{Th1.1}(i), we see that the tuple $(12,8,4,2,3,1)$ is universal.

From \eqref{Athm20u1}, we obtain
  \begin{equation}\label{Athm21u3}
     X\left(q^4\right)Y(q)Y\left(q^2\right)= X\left(q^4\right)X\left(q^8\right)Y\left(q^2\right)+q X\left(q^4\right)Y\left(q^2\right)Y\left(q^4\right).
  \end{equation}
  Since $(24,12,6,4,6,2)$ is universal by Sun \cite{Sun1}, with the aid of
  \eqref{2p58} we see that the tuple $(12,8,12,4,6,4)$ is also universal. Combining this with the second term of the identity  \eqref{Athm21u3} and Theorem \ref{Th1.1}(i), we obtain the universality of the tuple $(12,8,6,4,6,2)$. Hence $(16,8,8,4,6,4)$ is also universal in light of \eqref{SunLemma2}.

In view of \eqref{Athm20u1}, we have
\begin{equation}\label{Athm21u5}
  X\left(q^2\right)X\left(q^6\right)Y\left(q\right)= X\left(q^2\right)X\left(q^6\right)X\left(q^8\right)  +qX\left(q^2\right)X\left(q^6\right)Y\left(q^4\right).
\end{equation}
Since  the tuple $(18,6,6,4,6,2)$  is universal by Theorem \ref{Th2.4},
  applying Theorem \ref{Th1.1}(i) to the second term of \eqref{Athm21u5} we get the universality of the tuple $(12,8,9,3,3,1)$.

To determine the universality of the tuple $(12,8,12,6,3,1)$, using \eqref{Athm20u1} we have
  \begin{equation}\label{Athm20u5}
 \psi(q^6)X(q^2) Y(q)=\psi(q^6)X(q^2)X(q^8)+q\psi(q^6)X(q^2)Y(q^4).
  \end{equation}
  Sun \cite{Sun1} proved that the sum $x(3x+1)+y(3y+2)+z(3z+3)$ (or $(24,12,6,4,6,2)$) is universal. Thus, in light of
  the second term of  \eqref{Athm20u5} and Theorem \ref{Th1.1}(i), we obtain the universality of
$(12,8,12,6,3,1)$.

Using \eqref{appli12}, we have
\begin{equation}\label{univ2}
  \varphi(q)\psi\left(q^{4}\right) X\left(q^{2}\right)=\varphi\left(q^4\right)\psi\left(q^{4}\right) X\left(q^{2}\right)+2q\psi\left(q^{4}\right) \psi\left(q^8\right)X\left(q^{2}\right).
\end{equation}
By the paragraph containing \eqref{Athm23u6}, the tuples $(8,4,4,0,3,1)$ and $(16,8,8,4,3,1)$ are both universal. Thus, applying Theorem \ref{Th1.1}(i) to the identity
 \eqref{univ2} we get the universality of the tuple $(16,8,6,2,2,0)$.

 In view of the above, we have completed the proof of Theorem \ref{Th1.2}. \qed
\medskip

As mentioned in Section 2, the tuples $(20,10,5,1,3,1)$ and $(20,10,5,3,3,1)$
were proved to be universal in \cite[Theorem 1.1(iv)]{Sun3} via the theory of ternary quadratic forms.
Now we give a new proof by using our approach via theta functions.
 Using \cite[Corollary 2.11]{Cao} with $a=q$, $b=q^4$, $c=q^2$ and $d=q^3$, we obtain
\begin{align*}
f(q,q^4) f(q^2,q^3)= &f^2(q^{10},q^{15}) +q f(q^{15},q^{10}) f(q^{20},q^{5})+q^2 f^2(q^{20},q^{5}) \nonumber \\  &+2q^3 \psi (q^{25} ) f(q^{15},q^{10})+2q^4 \psi (q^{25})  f(q^{20},q^{5}),
\end{align*}
and hence
\begin{align}\label{Athm25u1}
X(q^{5}) f(q,q^4) f(q^2,q^3)= &X(q^{5}) f^2(q^{10},q^{15}) +qX(q^{5})  f(q^{15},q^{10}) f(q^{20},q^{5})\nonumber \\  &+q^2 X(q^{5}) f^2(q^{20},q^{5})   +2q^3\psi (q^{5})  X(q^{25} ) f(q^{15},q^{10})\nonumber \\  &+2q^4 X(q^{5}) \psi (q^{25})  f(q^{20},q^{5}),
\end{align}
Wu and Sun \cite{WuSun} showed that the tuple $(15,5,4,2,4,2)$ is universal.
Thus the tuple  $ (15, 5, 5, 3, 5, 1)$ is also universal with the aid of \eqref{513}.
 Therefore, in light of the last two  terms of the identity \eqref{Athm25u1} and Theorem \ref{Th1.1}(i), we see that the tuples $(20,10,5,1,3,1)$ and $(20,10,5,3,3,1)$ are universal.

\section{Proof of Theorem \ref{Th1.4}} \label{AlmostSec}

\medskip
\noindent{\tt Proof of Theorem \ref{Th1.4}}.
(i) Kane and Sun \cite{KaneSun} showed that the sum
 $$x^2+\frac{y(y+1)}{2}+9\frac{z(z+1)}{2}$$
 is almost universal. Note that
   $$x^2+\frac{y(y+1)}{2}+9\frac{z(z+1)}{2} \sim 9 x(2x+1)+y(3y+1)+\frac{z(3z+1)}{2} $$
   by \eqref{SunLemma1} and \eqref{p6}. So the tuple $(36,18,6,2,3,1)$ is almost universal.

    Let $|q|<1$. Multiplying both sides of \eqref{aAthm13u1}
by $\psi(q^9)=f(q^9,q^{27})$, we find that
\begin{equation}\label{almThm2.1}\begin{aligned}&\ f(q^9,q^{27})f(q,q^2)f(q^2,q^4)
\\=&\ f(q^9,q^9)f(q^9,q^{27})f(q^3,q^6)
+qf(q^9,q^{27})f(q^3,q^6)f(q^3,q^{15})
\\&\ +2q^2f(q^9,q^{27})^2f(q^3,q^{15}).
\end{aligned}
\end{equation}
As the tuple $(36,18,6,2,3,1)$ is almost universal, applying Theorem \ref{Th1.1}(ii)
to the identity \eqref{almThm2.1} and recalling \eqref{p6}, we find that
 the sums
$$3x^2+3\f{y(y+1)}2+\f{z(3z+1)}2\, \quad  \text{and}\,\quad  3\f{x(x+1)}2+\f{y(3y+1)}2+z(3z+2),$$
 are almost universal in view of the first and second terms of the right-hand side of \eqref{almThm2.1}.

 (ii) Kane and Sun \cite{KaneSun} showed that the sums
 $$x^2+6y(y+1)+\frac{z(z+1)}{2}\ \ \t{and}\ \ x^2+6y^2+\frac{z(z+1)}{2}$$
are almost universal.  By \eqref{SunLemma1} and \eqref{p6},
  $$x^2+6y(y+1)+\frac{z(z+1)}{2} \sim 12x(2x+1)+y(3y+1)+\frac{z(3z+1)}{2}$$
  and
  $$x^2+6y^2+\frac{z(z+1)}{2} \sim 6x^2+y(3y+1)+\frac{z(3z+1)}{2}.$$
   So the tuples $(24,12,6,2,3,1)$
  and $(12,0,6,2,3,1)$ are almost universal.

Let $|q|<1$. Multiplying both sides of \eqref{aAthm13u1} by $\psi(q^{12})=f(q^{12},q^{36})$, we obtain the identity
\begin{equation}\label{almThm2.3}\begin{aligned}
 &\ f(q^{12},q^{36})f(q,q^2)f(q^2,q^4)
 \\=&\ f(q^9,q^9)f(q^{12},q^{36})f(q^{3},q^{6})+qf(q^{12},q^{36})f(q^3,q^6)f(q^3,q^{15})
 \\&\ +2q^2f(q^9,q^{27})f(q^{12},q^{36})f(q^3,q^{15}).
\end{aligned}
\end{equation}
As the tuple $(24,12,6,2,3,1)$ is almost universal, applying Theorem \ref{Th1.1}(ii)
to the identity \ref{almThm2.3} and recalling \eqref{p6}, we see that the sums
   $$3x^2+2y(y+1)+\f{z(3z+1)}2\ \quad  \text{and}\,\quad 2x(x+1)+\f{y(3y+1)}2+z(3z+2),$$
    are almost universal in light of the first and second terms of the right-hand side of \eqref{almThm2.3}.
Similarly, as the tuple $(12,0,6,2,3,1)$ is almost universal. applying Theorem \ref{Th1.1}(ii)
to the identity \eqref{aAthm13u1} multiplied by $\varphi(q^6)=f(q^6,q^6)$, we deduce that
the sums
     $$2x^2+3y^2+\f{z(3z+1)}2 \quad  \text{and}\,\quad 2x^2+\f{y(3y+1)}2+z(3z+2),$$
       are both almost universal.

(iii) Kane and Sun \cite{KaneSun} proved that the sums
  $$3x(x+1)+y(y+1)+\f{z(z+1)}2\ \ \t{and}\ \ 6x^2+y(y+1)+\f{z(z+1)}2 $$
  are almost universal.
   By \eqref{SunLemma2} and \eqref{p6}, we have
   $$3x(x+1)+y(y+1)+\f{z(z+1)}2  \sim 6x(2x+1)+\f{y(3y+1)}2+{z(3z+2)} $$
   and
   $$ 6x^2+y(y+1)+\f{z(z+1)}2  \sim 6x^2+\f{y(3y+1)}2+{z(3z+2)}. $$
   Thus the tuples $(24,12,6,4,3,1)$ and $(12,0,6,4,3,1)$ are almost universal.

 Let $|q|<1$. Multiplying both sides of \eqref{aAthm13u6} by $\psi(q^6)=f(q^6,q^{18})$, we get the identity
 \begin{equation} \label{almThm2.4}
 \begin{aligned}&\ f(q^6,q^{18})f(q,q^2)f(q,q^5)
 \\=&\ f(q^6,q^{18})f(q^3,q^6)f(q^6,q^{12})+2qf(q^6,q^{18})f(q^9,q^{27})f(q^6,q^{12})
 \\&\ +2q^2f(q^6,q^{18})f(q^{18},q^{54})f(q^3,q^6).
 \end{aligned}
 \end{equation}
 As the tuple $(24,12,6,4,3,1)$ is universal, applying Theorem \ref{Th1.1}(ii) to the identity
 \eqref{almThm2.4} and recalling \eqref{p6}, we find that
  the sums
  $$x(x+1)+y(3y+1)+\f{z(3z+1)}2 \quad  \text{and}\,\quad x(x+1)+3y(y+1)+\f{z(3z+1)}2$$
  are almost universal in view of the first and third terms of the right-hand side of \eqref{almThm2.4}. Similarly, as the tuple $(12,0,6,4,3,1)$ is almost universal,
   applying Theorem \ref{Th1.1}(ii) to the identity \eqref{aAthm13u6} multiplied by
   $\varphi(q^6)=f(q^6,q^6)$ we deduce that the sums
 $$2x^2+y(3y+1)+\f{z(3z+1)}2 \quad  \text{and}\,\quad 2x^2+3y(y+1)+\f{z(3z+1)}2$$
  are both almost universal.

  In view of the above, we have completed the proof of Theorem \ref{Th1.4}. \qed

\end{document}